\documentclass[12pt]{article}
\usepackage{amsmath}
\usepackage{amssymb}
\usepackage{amsthm}
\usepackage{mathtools}
\usepackage{dsfont}
\usepackage{graphicx} 
\usepackage{tikz}
\usetikzlibrary{shapes}
\usepackage{faktor}
\usepackage{authblk}
\usepackage{enumerate}
\usepackage{color}
\usepackage{overpic}
\usepackage[utf8]{inputenc}
\usepackage[
    inner=2.3cm,
    outer=2.3cm,
    top=2.7cm,
    marginparwidth=2.5cm,
    marginparsep=0.1cm
    ]
    {geometry}
    \usepackage[hidelinks]{hyperref}

\newcommand{\X}{\mathds{X}}
\newcommand{\Z}{\mathds{Z}}
\newcommand{\R}{\mathds{R}}
\newcommand{\N}{\mathds{N}}
\newcommand{\E}{\mathds{E}}

\newcommand{\abs}[1]{\left| #1 \right|}
\newcommand{\1}{\mathds{1}}
\newcommand{\supp}{\text{supp}}
\newcommand{\F}{\mathcal{F}}
\renewcommand{\P}{\mathds{P}}
\newenvironment{nalign}{
    \begin{equation}
    \begin{aligned}
}{
    \end{aligned}
    \end{equation}
    \ignorespacesafterend
}

\newtheorem{theorem}{Theorem}[section]
\newtheorem{proposition}[theorem]{Proposition}

\newtheorem{corollary}[theorem]{Corollary}
\newtheorem{definition}[theorem]{Definition}
\newtheorem{conjecture}{Conjecture}
\newtheorem*{theorem*}{Theorem}

\theoremstyle{definition}
\newtheorem{example}{Example}
\newtheorem{remark}[theorem]{Remark}

\title{Random attractors on countable state spaces}
\author[1]{Robin Chemnitz}
\author[1,2]{Maximilian Engel}
\author[3]{Guillermo Olicón Mendez}

\affil[1]{Freie Universität Berlin}
\affil[2]{University of Amsterdam}
\affil[3]{Universität Potsdam}

\date{May 2024}

\begin{document}
\maketitle

\begin{abstract}
We study the synchronization behavior of discrete-time Markov chains on countable state spaces. Representing a Markov chain in terms of a random dynamical system, which describes the collective dynamics of trajectories driven by the same noise, allows for the characterization of synchronization via random attractors.   

We establish the existence and uniqueness of a random attractor under mild conditions and show that forward and pullback attraction are equivalent in our setting.
Additionally, we provide a sufficient condition for reaching the random attractor, or synchronization respectively, in a time of finite mean.

By introducing insulated and synchronizing sets, we structure the state space with respect to the synchronization behavior and characterize the size of the random attractor. 
\end{abstract}

\section{Introduction}
In this work, we consider discrete-time Markov chains on a countable state space which are irreducible and positive recurrent, and, hence, admit a unique stationary distribution. The law of a Markov chain $(X_n)_{n\in\N}$ is governed by the transition probabilities $\P(X_{n+1} = y \mid X_n=x)$. These probabilities describe the random motion of a single particle, but contain no information about how multiple particles move collectively, i.e.~how the individual trajectories are coupled. We consider a specific type of coupling, which arises from representing a Markov chain by a \emph{random dynamical system}, short RDS, \cite{Arnold1998}. 

Random dynamical systems model in a pathwise manner the evolution of stochastic processes for fixed noise realizations, allowing the simultaneous tracking of multiple trajectories under the same randomness. RDS representations of Markov chains on discrete state spaces have previously been studied in \cite{ye2016, grune2021random, EngelOlicon24}. Informally, an RDS representation of a Markov chain with state space $\X$ randomly picks a map $f_n:\X\to \X$ in each time step $n\in \N$ such that the process $X_{n+1} = f_n(X_n)$ admits the law of the given Markov chain. For a formal definition of RDS representations of a Markov chain, we refer to Section \ref{SEC:preliminaries}. Note that such a representation naturally describes the collective behavior of any number of particles, i.e.~$Y_{n+1} = f_n(Y_n)$ for another starting condition $Y_0$. Note that if $X_n=Y_n$ for some $n\in \N$, then $X_m = Y_m$ for all $m\geq n$. This is what we call \emph{synchronization}. If any two initial conditions $X_0$, $Y_0$ eventually synchronize, we call the RDS \emph{synchronizing}. We formalize the synchronization behavior of trajectories using the notion of a \emph{random attractor}. 

\subsubsection*{Related work}

For RDS on general state spaces, the theory of attractors has been thoroughly developed over the past decades, see for instance \cite{CrauelFlandoli94, crauelFlandoli97, Arnold1998, crauel2015, Crauel16}. 
Generally, different notions of attraction coexist without redundancy, where examples have been typically studied in the context of stochastic differential equations \cite{Scheutzow2002, ChueshovScheutzow2004, CrauelDimitroffScheutzow2009, CrauelGeissScheutzow2024, DimitroffScheutzow2011}.
In more detail: as opposed to the autonomous deterministic setting,  random dynamical systems (being non-autonomous) may exhibit different notions of attraction, e.g.~\textit{pullback} and \textit{forward attractors} \cite{RasmussenKloeden11}. Additionally, random attractors are random sets whose notion of attraction is further classified by the type of convergence, i.e.~in the strong or the weak sense \cite{crauel2015}. One can also distinguish random attractors  in terms of the family of subsets which are attracted towards them.
We show that in our setting, which considers discrete state space, almost all of the commonly used notions of attraction are equivalent, cf.~Theorem~\ref{thm:equivalence} and Remark \ref{rem:equivalence} thereafter. This is a direct generalization of results formulated in \cite{EngelOlicon24}.

The representation of a Markov chain by an RDS is very similar to the notion of a \emph{stochastic recursive sequence} \cite{borovkov1992stochastically}. More recently, the representation of Markov chains by stochastic recursive sequences has been studied using the framework of unimodular networks \cite{Baccellietal2018, Baccellietal2019}. One of the key motivations for these works is the \emph{coupling from the past} algorithm, which was introduced in \cite{propp1996exact} and provides a method to draw perfect samples from the stationary distribution of a Markov chain. The objects we consider in this work are similar to those studied in the references mentioned above, despite being formulated in a different mathematical language. For example, a random attractor corresponds to a finite set of bi-infinite paths in \cite{Baccellietal2019}. Using the framework of random dynamical systems, which enables the use of powerful tools like the correspondence theorem, cf.~Theorem \ref{thm:correspondence}, we obtain the existence of a unique random attractor, characterize its size, and derive a detailed structuring of the state space into insulated and synchronizing sets, cf.~Section \ref{sec:structure_of_attractor}. 
Furthermore, Section \ref{sec:time_until_attraction} concerns the expected time until attraction, or synchronization respectively. This quantity has previously been considered in \cite{Baccellietal2019}. We provide a sufficient condition for the expected hitting times to be finite.

\subsubsection*{Main results and structure of the paper}

The main results of our work are the following.
\begin{itemize}
    \item Theorem \ref{thm:equivalence}. There exists a unique random attractor that is both a pullback attractor for finite deterministic sets and a forward attractor for compact random sets. This shows that for finite deterministic sets forward and pullback attraction are equivalent.
    \item Theorem \ref{thm:synchronization_finite_mean}. For synchronizing RDS, the expected time until two given initial conditions synchronize is always finite.
    \item Theorem \ref{thm:finite_mean_iff_erg_2}. If the associated Markov chain is ergodic of degree 2, the expected time to hit the attractor is finite for every initial condition.
    
\end{itemize}

In Section \ref{SEC:preliminaries}, we recall some theory of Markov chains on countable state spaces and formally introduce random dynamical systems and different notions of random attractors. Section \ref{sec:equivalence} is dedicated to the proof of Theorem \ref{thm:equivalence} and establishes the existence of a unique random attractor $A$ that has a fixed size $\kappa$. In Section \ref{sec:structure_of_attractor} we structure the state space with respect to the synchronizing behavior and introduce insulated and synchronizing sets. In particular, we show that the random attractor $A$ is a maximum insulated set (Theorem~\ref{thm:A_anti_clique}), and derive restrictions to the possible values of $\kappa$ (cf.~Corollary~\ref{crl:restriction_kappa}). In Section \ref{sec:time_until_attraction}, we study how quickly initial conditions hit the attractor, as well as how quickly synchronization occurs in the case of synchronizing RDS. In particular, we prove Theorem~\ref{thm:finite_mean_iff_erg_2} and Theorem~\ref{thm:synchronization_finite_mean}, which provide finite expected times until attraction, or synchronization respectively. At the end of each section, we provide several examples illustrating the results of the respective section. Finally, we discuss our findings in Section \ref{sec:discussion} where we also address some open questions.

\section{Preliminaries}
\label{SEC:preliminaries}
\subsection{Markov chains}
We consider time-homogeneous discrete-time Markov chains $(X_n)_{n\in \N}$ on a countable state space $\X$, where $\N$ includes $0$. A Markov chain is defined by its transition probabilities for each pair $x,y\in \X$
\begin{equation}\label{eq:trans_prob}
    P(x,y):= \P(X_1=y \:|\: X_0=x), \quad P^n(x,y):= \P(X_{n}=y \:|\: X_0=x).
\end{equation}
A stationary distribution of $(X_n)_{n\in \N}$ is a probability vector $\pi\in \R^\X_{\geq 0}$ that satisfies
\begin{equation*}
    (\pi P) (y) := \sum_{x\in \X}\pi(x) P(x,y) = \pi(y), \quad \forall y\in \X, \quad \text{and} \quad \sum_{x\in \X} \pi(x)=1.
\end{equation*}
We call a Markov chain \textit{irreducible} if for each $x,y\in \X$ there is an $n \in \N$ such that $P^n(x,y)>0$. An \emph{irreducibility class} $I\subset \X$ is a set such that the Markov chain restricted to $I$ is irreducible, i.e.~for all $x,y\in I$ there is an $n\in \N$ such that $P^n(x,y)>0$, and for all $x\in I$ and $y\in \X \setminus I$ we find $P^n(x,y)=0$ for all $n\in \N$.

Denote the first hitting time of a state $x\in \X$ by  
\begin{equation}\label{eq:MC_hitting_time}
    \tau_x := \min \{n\geq 1 \: | \: X_n = x\}.
\end{equation}
We denote the probability and the expectation, conditioned on $X_0=x$, by $\P_x$ and $\E_x$ respectively.
Hence, the expected return time of a state $x\in\X$ is written as $\E_x[\tau_x]$. We call the Markov chain \emph{recurrent} if
\begin{equation*}
    \P_x(\tau_x < \infty) = 1, \quad \forall x\in \X.
\end{equation*}
We call the Markov chain \emph{positive recurrent} if the expected return times are of finite mean for all states, i.e.
\begin{equation*}
    \E_x[\tau_x] < \infty, \quad \forall x\in \X.
\end{equation*}

The following equivalence is well-known, see e.g.~\cite[Theorem 2.6]{benaim2022markov}.
\begin{theorem}\label{thm:MC_ergodic}
    For an irreducible Markov chain $(X_n)_{n\in \N}$ on a countable state space, the following notions are equivalent:
    \begin{enumerate}[(i)]
        \item There is an $x\in \X$ such that $\E_x[\tau_x] < \infty$; 
        \item $(X_n)_{n\in \N}$ is positive recurrent;
        \item $(X_n)_{n\in \N}$ admits a unique stationary distribution $\pi$ that satisfies $\pi(x) = \E_x[\tau_x]^{-1}>0$, for any $x\in \X$.
    \end{enumerate}
    Under these equivalent conditions:
    \begin{itemize}
        \item For any $x,y\in \X$ we find the $\P_x$-a.s.~limit
        \begin{equation}\label{eq:Markov_chain_Birkhoff}
            \lim_{N\to \infty} \frac{1}{N}\sum_{n=1}^N \1_y(X_n) = \pi(y);
        \end{equation}
        \item $\E_x[\tau_y] < \infty$, for all $x,y\in \X$.
    \end{itemize}
\end{theorem}
In the following, we always work with irreducible and positive recurrent Markov chains whose unique stationary distribution is denoted by $\pi$. Note that we do not assume aperiodicity.

\subsection{Random dynamical systems}
We briefly introduce the general notion of a random dynamical system (RDS), before focusing on the representation of a Markov chain by an RDS. Random dynamical systems model the evolution of stochastic processes pathwise for fixed noise realizations. This allows for simultaneously tracking multiple trajectories under the same randomness, giving rise to interesting phenomena like random attractors. 
Since we only consider countable state spaces, we directly provide all definitions in the setting of a countable space $\X$ which we equip with the powerset $\sigma$-algebra, denoted by $\mathcal{P}(\X)$. 
For a detailed introduction to random dynamical systems on general state spaces in discrete and continuous time, we refer to \cite{Arnold1998}.

An RDS $(\theta, \varphi)$ on $\X$ consists of a probability space $(\Omega, \F, \P)$, a group of bimeasurable measure-preserving transformations $(\theta_n:\Omega \to \Omega)_{n\in \Z}$ such that for all $\omega \in \Omega$
\begin{align}
\begin{split}
     \theta_0 \omega &= \omega,\\
     \theta_{n+m} \omega &= (\theta_n \circ \theta_m) \omega, \quad \forall m,n \in \Z,
 \end{split}
 \end{align}
 and a $(\mathcal{F}\otimes \mathcal{P}(\mathbb{N}) \otimes \mathcal{P}(\X))$-measurable map $\varphi: \Omega \times \N \times \X \to \X$ that has the cocycle property
 \begin{align}
     \begin{split}
         \varphi_\omega^0 (x) &= x,\\
         \varphi^{n+m}_\omega (x) &= \left( \varphi^n_{\theta_m \omega} \circ \varphi^m_\omega \right) (x),
     \end{split}
 \end{align}
 for all $\omega \in \Omega$, $m,n\in \N$ and $x\in \X$. Due to the cocycle property, an RDS is uniquely defined by the time-one maps $\varphi_\omega^1:\X \to \X$ and $\theta_1:\Omega \to \Omega$. We usually write $\varphi_\omega$ and $\theta$ instead of $\varphi_\omega^1$ and $\theta_1$. We additionally assume that $\theta$ is ergodic with respect to $\P$, which implies that $\theta_{-1}$ is ergodic with respect to $\P$ as well. From the cocycle property it follows that if two initial states $x,y\in \X$ coalesce, i.e.~$\varphi_\omega^n(x) = \varphi_\omega^n(y)$ for some $n\in \N$, then $\varphi_\omega^m(x) = \varphi_\omega^m(y)$ for all $m>n$. We say that $x$ and $y$ $\omega$-\emph{synchronize} \cite{huang2019}.

The maps $\varphi_\omega^n$ can also be applied to sets and tuples. For $K$ being a set or tuple of elements in $\X$, we write $\varphi_\omega^n(K)$ for the set or tuple consisting of $\varphi_\omega^n(x)$, $x\in K$. Whether $\varphi^n_\omega$ is applied to an element, a set or a tuple should be apparent from the context.

Next, we define invariant measures of an RDS. 
The RDS $(\theta, \varphi)$ induces the so-called skew-product map on $\Omega \times \X$ given by
\begin{align*}
    \Phi: \Omega \times \X &\to \Omega \times \X, \\
    (\omega, x) &\mapsto (\theta \omega, \varphi_\omega(x)).
\end{align*}
An invariant measure of the RDS $(\theta, \varphi)$ will be defined to be invariant with respect to this skew-product map. For any measure $\mu$ on a measurable space $M$ and a measurable map $f:M\to M$, we define the push-forward of $\mu$ by $f$,
\begin{equation*}
    f_* \mu (\cdot) = \mu(f^{-1}(\cdot) ).
\end{equation*}
Hence, the push-forward operator $f_*$ maps measures on $M$ to measures on $M$. In particular, $f_*$ maps probability measures to probability measures. 
\begin{definition}\label{def:invariant_measure}
    A probability measure $\mu$ on $(\Omega \times \X, \F \otimes \mathcal{P}(\X))$ is called \emph{invariant} if
    \begin{enumerate}[(i)]
        \item $\Phi_* \mu = \mu$;
        \item The marginal of $\mu$ on $\Omega$ is $\P$ such that $\mu$ can be factorized uniquely into
        \begin{equation*}
            \mu(d\omega, x) = \mu_\omega(x) \P(d\omega),
        \end{equation*}
        where the $\mu_\omega$ are the sample measures on $\X$, i.e.
        \begin{itemize}
            \item $\mu_\omega$ is a probability measure on $\X$ for $\P$-a.e.~$\omega \in \Omega$;
            \item $\omega \mapsto \mu_\omega(x)$ is measurable for any $x\in \X$. 
        \end{itemize}
    \end{enumerate}
\end{definition}
Invariance of a measure $\mu$ on $\Omega \times \X$ with marginals $\P$ can be checked directly via the invariance of the sample measures. This alternative characterization tends to be much easier to work with.
\begin{proposition}
    A probability measure $\mu$ on $(\Omega \times \X, \F \otimes \mathcal{P}(\X))$ with marginal $\P$ on $\Omega$ is invariant with respect to the RDS $(\theta, \varphi)$, if and only if the samples measures satisfy
    \begin{equation}\label{eq:def_invariant_sample_measures}
        (\varphi_\omega)_* \mu_\omega = \mu_{\theta \omega},
    \end{equation}
    for $\P$-a.e.~$\omega\in \Omega$.
\end{proposition}
A proof of this result can be found in \cite[Proposition 1.3.27]{KuksinShirikyan2012}. Note that due to the cocycle property, equation \eqref{eq:def_invariant_sample_measures} implies that $(\varphi^n_\omega)_* \mu_\omega = \mu_{\theta_n \omega}$ for all $n\in \N$ and $\P$-a.e.~$\omega \in \Omega$.

\subsection{Markov RDS}\label{sec:Markov_RDS}
Let $(\theta, \varphi)$ be an RDS on a countable state space $\X$. For any $s,t\in \Z$ with $s<t$, we denote by $\F_s^t\subset \F$ the $\sigma$-algebra generated by (the subsets of $\F$ of zero measure and) the random variables $\omega \mapsto \varphi^n_{\theta_m \omega}(x)$ for $x\in \X$ and $m,n\in \Z$ with $s \leq m \leq m+n \leq t$. We define the $\sigma$-algebras
\begin{nalign}\label{eq:def_F_infty}
    \F_{-\infty}^t &:= \sigma (\F_s^t\: | \: s \in \Z, s<t), \\ 
    \F_s^{\infty} &:= \sigma (\F_s^t\: | \: t \in \Z, s<t), \\
    \F_{-\infty}^\infty &:= \sigma(\F_{-t}^t \: | \ \ t\in \N_{>0}).
\end{nalign}

\begin{definition}
    An RDS $(\theta, \varphi)$ is called \emph{Markov} if $\F_{-\infty}^0$ is independent of $\F_0^\infty$.
\end{definition} 

The $\sigma$-algebra $\F_{-\infty}^\infty$ contains all the dynamical information of the RDS. Clearly, $\F_{-\infty}^\infty \subset \F$, but in general $\F$ might contain strictly more sets than $\F_{-\infty}^\infty$. 

Note that for a Markov RDS, $\F_{-\infty}^t$ is independent of $\F_t^\infty$ for any $t\in \Z$. In particular, the $\sigma$-algebras $\F_t^{t+1}$ are mutually independent, and the maps $\varphi_\omega$ and $\varphi_{\theta_n \omega}$ are independent for all $n\in\Z \setminus \{0\}$. These observations imply that, given a Markov RDS $(\theta, \varphi)$, the random process defined by
\begin{equation}\label{eq:RDS_MC}
    X_0(\omega) = x, \quad X_n(\omega) = \varphi^n_\omega( x),
\end{equation}
is a time-homogeneous Markov chain \cite[Theorem 2.1.4]{Arnold1998}. The transition probabilities $P$ of $X$ are given by
\begin{equation*}
    P(x, y) = \P(\varphi_\omega(x) = y).
\end{equation*}
We say that the RDS represents the Markov chain $(X_n)_{n\in \N}$. Conversely, we call $(X_n)_{n\in \N}$ the associated Markov chain of the RDS. The key quality of an RDS representation is that it models the simultaneous motion of arbitrarily many points in $\X$ under the same noise realisation $\omega \in \Omega$, while the Markov chain only describes the random motion of a single point in $\X$.
\begin{remark}
    While each Markov RDS represents a unique Markov chain, a given Markov chain with transition probabilities $P$ has, in general, many possible RDS representations, see e.g.~\cite{ye2016}. Most of our results, like existence and uniqueness of a random attractor, cf.~Theorem \ref{thm:equivalence}, and the finiteness of expected hitting times studied in Section \ref{sec:time_until_attraction}, do not depend on the chosen RDS representation. However, quantities like the $\P$-a.s.~size of the attractor do depend on the RDS representation, cf.~Example \ref{ex:4_state_MC}.
\end{remark}

For Markov RDS, we strengthen the notion of an invariant measure, cf.~Definition \ref{def:invariant_measure}, and introduce the notion of an invariant Markov measure.
\begin{definition}
    Let $(\theta ,\varphi)$ be a Markov RDS. An invariant measure $\mu$ is called an \emph{invariant Markov measure} if the sample measures $\mu_\omega$ are $\F_{-\infty}^0$-measurable, i.e.~for each $x\in \X$, the map $\omega \mapsto \mu_\omega(x)$ is $\F_{-\infty}^0$-measurable.
\end{definition}

A central tool for our work is the \emph{correspondence theorem}, which links invariant Markov measures of a Markov RDS to the stationary distributions of the associated Markov chain in a one-to-one correspondence. In particular, if the associated Markov chain has a unique stationary distribution, the representing RDS has a unique invariant Markov measure. For a proof of the following theorem, see \cite[Theorem 4.2.9]{KuksinShirikyan2012}.
\begin{theorem}[Correspondence Theorem]\label{thm:correspondence}
    Let $(\theta, \varphi)$ be a Markov RDS. The following assertions hold:
    \begin{enumerate}[(i)]
        \item Let $\mu$ be an invariant Markov measure on $\Omega \times \X$ with sample measures $\mu_\omega$. Then 
        $$\pi:=\E[\mu_\omega]$$
         is a stationary distribution of the associated Markov chain. \\
        Additionally, if $\mu'$ is another invariant Markov measure with sample measures $\mu'_{\omega}$ such that $\E[\mu_\omega]=\E[\mu'_\omega]$, then $\mu=\mu'$.
    \item Let $\pi$ be a stationary distribution of the associated Markov chain. Then the weak limit
        \begin{equation}\label{eq:correspondence_limit}
            \mu_\omega := \lim_{n\to \infty} (\varphi^{n}_{\theta_{-n}\omega} )_* \pi,
        \end{equation}
        exists for $\P$-a.e.~$\omega \in \Omega$ and the measure $\mu$ defined by the sample measures $\mu_\omega$ is an invariant Markov measure such that $\E[\mu_\omega]= \pi$. 
    \end{enumerate}
\end{theorem} 

Note that the correspondence theorem only establishes a one-to-one correspondence between stationary distributions of the associated Markov chain and invariant \emph{Markov} measures on $\Omega \times \X$. In general, there can be other, non-Markovian invariant measures on $\Omega \times \X$. Under the assumption $\F = \F_{-\infty}^\infty$, \cite[Theorem 1.7.2]{Arnold1998} shows that the invariant Markov measure given by $\mu_\omega = \lim_{n\to \infty} (\varphi^{n}_{\theta_{-n}\omega} )_* \pi$ is ergodic with respect to the skew-product map $\Phi$ if and only if $\pi$ is an ergodic stationary distribution of the associated Markov chain. In particular, an RDS representing an irreducible and positive recurrent Markov chain has a single invariant measure. This measure is Markovian and ergodic with respect to $\Phi$.

\subsection{Prototypical constructions}\label{sec:canonical_constructions}
We introduce well-known constructions for Markov RDS representing a given Markov chain. We will use these constructions in most of the examples in the following sections.

Firstly, we introduce the \emph{shift space}.
Given a Markov chain on a countable state space $\X$ with transition probabilities $P$, consider the set $\Sigma$ of all maps $f:\X \to \X$. We equip $\Sigma$ with the $\sigma$-algebra $\mathcal{G}$ generated by the maps $f\mapsto f(x)$, for $x\in \X$. Let $\P_0$ be a probability measure on $(\Sigma, \mathcal{G})$ satisfying $\P_0(f(x) = y) = P(x,y)$ (see, for instance, equation~\eqref{eq:construction} below).
Define the probability space $(\Omega, \F, \P):=(\Sigma^\Z, \mathcal{G}^{\otimes \Z}, \P_0^{\otimes \Z})$. We write elements of $\Omega$ as $\omega=(\hdots, \omega_{-1}, \omega_0, \omega_1, \hdots)$, where $\omega_n \in \Sigma$.  
It is easy to observe that the left-shift $\theta:\Omega \to \Omega$, i.e.~$(\theta\omega)_n = \omega_{n+1}$, preserves the measure $\mathbb{P}$. 
Over this shift space $(\Omega, \F, \P, \theta)$, we define a cocycle by $\varphi_\omega(x) = \omega_0(x)$ and observe that  the resulting RDS $(\theta, \varphi)$ is a Markov RDS representation of the Markov chain on $\X$ with transition probabilities $P$. 
The measure $\P_0$ on $\Sigma$ uniquely determines the behavior of the RDS representation.

A special case of the construction above is the \emph{independent RDS}. In the independent RDS points move independently from one to another until they coalesce. Using Carath\'eodory's extension theorem \cite[Theorem 1.11]{Folland99}, define a measure $\P_0$ on $(\Sigma, \mathcal{G})$ by
\begin{equation} \label{eq:construction}
    \P_0(f(K) = L) = \prod_{i=1}^n P(x_i, y_i),
\end{equation}
where $K=(x_1, \hdots, x_n)$ and $L=(y_1, \hdots, y_n)$ are tuples of distinct elements of $\X$. We call the resulting RDS over the shift-space $(\Sigma^\Z, \mathcal{G}^{\otimes \Z}, \P_0^{\otimes \Z})$ the independent RDS. We will present and analyze the independent RDS of a given Markov chain in Example \ref{ex:4_state_MC}. 

\subsection{Random attractors}
Let $(\theta, \varphi)$ be a Markov RDS on a countable state space $\X$. A central object of study in the theory of random dynamical systems is the notion of a random attractor. 
In the literature, there are multiple notions of attractors that are in general neither equivalent nor strictly stronger than one another, cf.~\cite{Scheutzow2002}. All of these definitions introduce attractors as invariant compact random sets. For a detailed discussion on the matter, we refer to \cite{crauel2015}. Focusing on discrete time and a countable state space, we give the definitions of such sets accordingly. We equip $\X$ with the discrete topology such that a subset of $\X$ is compact if and only if it is finite. Hence, in the following, we use compact and finite interchangeably.
\begin{definition}
    We call $A:\Omega\to \mathcal{P}(\X)$ a compact random set if
    \begin{enumerate}[(i)]
        \item $\omega\mapsto \1_{A(\omega)}(x)$ is measurable for each $x\in \X$;
        \item $A(\omega)$ is compact, i.e.~finite, for $\P$-a.e.~$\omega \in \Omega$.
    \end{enumerate}
    We say that $A$ is an invariant compact random  set if additionally
    \begin{enumerate}
        \item[(iii)] $A$ is $\varphi$-invariant, i.e.
        \begin{equation*}
            \varphi_\omega (A(\omega)) = A(\theta \omega), \quad \P\text{-a.s.}.
        \end{equation*}
    \end{enumerate}
\end{definition}

\begin{definition}\label{def:attractors}
    An invariant compact random set $A$ is called a
    \begin{enumerate}[(i)]
        \item pullback attractor for compact random sets, if for each compact random set $B$
        \begin{equation*}
            \P(\exists N\in \N \,, \forall n\geq N: \varphi^n_{\theta_{-n}\omega} (B(\theta_{-n}\omega)) \subset A(\omega))= 1.
        \end{equation*}
        \item forward attractor for compact random sets, if for each compact random set $B$
        \begin{equation*}
            \P(\exists N\in \N \,, \forall n\geq N: \varphi^n_\omega (B(\omega)) \subset A(\theta
            _n \omega))= 1.
        \end{equation*}
        \item weak attractor for compact random sets, if for each compact random set $B$
        \begin{nalign}\label{eq:def_weak_attractor}
            \lim_{n\to \infty} \P(\varphi^n_{\theta_{-n}\omega} (B(\theta_{-n} \omega)) \subset A(\omega))&= 1, \\
            \lim_{n\to \infty} \P(\varphi^n_\omega (B(\omega)) \subset A(\theta_n\omega))&=
            1.
        \end{nalign}
        \item pullback/forward/weak attractor for finite deterministic sets, if the respective definition holds for finite deterministic sets $B\subset \X$;
        \item pullback/forward/weak point attractor, if the respective definition holds for all singleton sets $B=\{x\}$, $x\in \X$.
    \end{enumerate}
\end{definition}
\begin{remark}
    The two criteria in the definition of weak attraction are equivalent by the measure-preserving property of $\theta$. Both of these criteria are listed to emphasize that a weak attractor is both a weak forward and a weak pullback attractor. 
\end{remark}
For general RDS, pullback attraction and forward attraction are neither weaker nor stronger than the other, cf.~\cite{Scheutzow2002}. On the other hand, it follows directly from the definitions that both pullback and forward attraction are stronger than the respective weak attraction. Our results of the next section, cf.~Theorem \ref{thm:equivalence} and Example \ref{ex:no_pullback}, show that in our discrete setting pullback attraction is stronger than forward attraction.

In many settings, attraction of compact random sets is too strong of a notion. Hence, it is common to restrict attraction of compact random sets to a class of 'nice' sets. In \cite[Definition~14.3]{RasmussenKloeden11}, a pullback attractor is only required to attract so-called \emph{tempered} random sets. In the next section, we show that any Markov RDS representing an irreducible and positive recurrent Markov chain admits a forward attractor for general compact random sets. Therefore, we keep the definitions at this high level of generality. 

Closely related to the concept of attraction, we introduce the notion of \emph{synchronization}.
\begin{definition}\label{def:omega_sync}
   We say that $x,y\in \X$ \emph{$\omega$-synchronize} if there is an $n\in \N$ such that $\varphi^n_\omega(x) = \varphi^n_\omega(y)$. An RDS $(\theta, \varphi)$ is called \emph{synchronizing}, if any $x,y\in \X$ $\omega$-synchronize for $\P$-a.e.~$\omega \in \Omega$. 
\end{definition}
For a synchronizing RDS, we define the synchronization time of two points $x,y\in \X$ as
\begin{equation}\label{eq:synchronization_time}
    T_\Delta(\omega, x,y) := \min\{ n\in \N \:  | \: \varphi_\omega^n(x) = \varphi_\omega^n(y)\}.
\end{equation}
This time is $\P$-a.s.~finite by definition. In Section~\ref{sec:time_until_attraction}, we will show the also the \emph{expected} synchronization time is always finite, cf.~Theorem \ref{thm:synchronization_finite_mean}.

If an RDS has a forward point-attractor $A$ that $\P$-a.s.~consists of a single point, i.e.~$\abs{A(\omega)}=1$ for $\P$-a.e.~$\omega \in \Omega$, then the RDS is synchronizing. For Markov RDS, representing an irreducible and positively recurrent Markov chain, the converse holds true as well, cf.~Remark \ref{rem:kappa_sync}. In general, the converse is not true since there are RDSs that are synchronizing but do not have an attractor.

\section{Random attractors for Markov RDS}
The main result of this section states that an RDS representing an irreducible and positive recurrent Markov chain possesses a unique random attractor that satisfies all the notions of attraction in Definition \ref{def:attractors} except pullback attraction for compact random sets. In particular, this shows that, apart from pullback attraction for compact random sets, all notions of attraction are equivalent in this setting. Therefore, in this setting, pullback attraction is a strictly stronger notion of attraction than forward attraction. In Example \ref{ex:no_pullback} we illustrate that without further assumptions, the unique attractor is not a pullback attractor for compact random sets. We emphasize that the results of this section hold for any RDS representation of a given Markov chain. 

\subsection{Existence and uniqueness}\label{sec:equivalence}
\begin{theorem}\label{thm:equivalence}
    Let $(\theta, \varphi)$ be a Markov RDS on a countable state space $\X$ representing a discrete-time Markov chain which is irreducible and positive recurrent. 
    \begin{enumerate}
        \item[(a)] There is an $\F_{-\infty}^0$-measurable invariant compact random set $A$ which is a pullback attractor for finite deterministic sets and a forward attractor for compact random sets.
        \item[(b)] A is the unique attractor in the sense that $A$ is (up to nullsets) the only invariant compact random set satisfying some notion of attraction from Definition \ref{def:attractors}.
        \item[(c)] There is a constant $\kappa\in \N_{>0}$ such that $\abs{A(\omega)}= \kappa$ for $\P$-a.e.~$\omega\in \Omega$.
        \item[(d)] The probability that a state $x\in \X$ lies in the attractor is given by
    \begin{equation}\label{eq:prob_x_in_A}
        \P(x\in A(\omega)) = \kappa \pi(x).
    \end{equation}
    \end{enumerate} 
\end{theorem}

The rest of this subsection is dedicated to proving this theorem. Instead of writing out one unified proof of $(a)$, we formulate multiple smaller propositions that together provide a proof. The structure is as follows:
\begin{enumerate}[(I)]
    \item There is a weak point attractor $A$ which is $\F^0_{-\infty}$-measurable;
    \item $A$ is a forward point attractor;
    \item $A$ is a pullback point attractor;
    \item Forward point attraction is equivalent to forward attraction of 
    compact random sets and pullback point attraction is equivalent to pullback attraction for finite deterministic sets.
    \item The statements $(b)$, $(c)$ and $(d)$ hold.
\end{enumerate}
\begin{remark}\label{rem:equivalence}
    Parts $(a)$ and $(b)$ of Theorem \ref{thm:equivalence} show that for irreducible and positive recurrent Markov chains all notions of attraction from Definition \ref{def:attractors}, apart from pullback attraction for compact random sets, are equivalent. If there is an attractor $A$ of any type, then by the assertions of the theorem, this is the unique attractor of the RDS and thereby satisfies all other notions of attraction apart from pullback attraction for compact random sets. In Example \ref{ex:no_pullback}, we illustrate that in general $A$ is indeed not a pullback attractor for compact random sets.
\end{remark}

Throughout this section, we consider a fixed RDS $(\theta, \varphi)$ representing a discrete-time Markov chain $(X_n)_{n\in \N}$ with transition probabilities $P$, that is irreducible and positive recurrent with stationary distribution $\pi$.
We complete step (I) using results from \cite{Scheutzow17}. The existence of a weak attractor under strictly stronger assumptions has been shown in \cite{EngelOlicon24}. The following proposition is a direct extension of this result.
\begin{proposition}[I]\label{prop:existence_weak_attractor}
    There is a weak point attractor $A$ which is $\F^0_{-\infty}$-measurable.
\end{proposition}
\begin{proof}
    This theorem uses results from \cite{Scheutzow17}. In their work, $\F = \F_{-\infty}^\infty$ is assumed. Note that $(\theta, \varphi)$ is a well-defined RDS over the probability space $(\Omega, \F_{-\infty}^\infty, \P)$. Since $\F_{-\infty}^\infty \subset \F$, any point attractor of the RDS with $\sigma$-algebra $\F_{-\infty}^\infty$ is also a point attractor of the RDS with $\sigma$-algebra $\F$. Hence, w.l.o.g.~we assume $\F= \F_{-\infty}^\infty$ whenever we use results from \cite{Scheutzow17} to show the existence of a point attractor. 
    
    The correspondence theorem, in particular (\ref{eq:correspondence_limit}), provides a unique invariant Markov measure $\mu$ for $(\theta, \varphi)$ on $\Omega\times \X$.
    By \cite[Lemma 2.19]{Scheutzow17},
    we know that  the sample measures $\mu_\omega$ either $\P$-a.s.~have no point mass or $\P$-a.s.~consist of finitely many atoms of the same mass.
    Since we work on a countable state space that does not support probability measures without point mass, we can find a random set $A$ and $\kappa \in \N$ such that for $\P$-a.e.~$\omega\in \Omega$
    \begin{equation}\label{eq:sample_measures_uniform}
        \mu_\omega = \frac{1}{\kappa} \sum_{a\in A(\omega)} \delta_a,
    \end{equation}
   where  $|A(\omega)|=\kappa$.
    In particular, we have for such $\omega$ that
    \begin{equation}\label{eq:A_supp_mu}
        A(\omega) = \supp(\mu_\omega).
    \end{equation}
    Hence, $A$ is a random set that is $\P$-a.s.~compact.
    Note that there is a $\theta$-invariant set $\hat{\Omega}\subset \Omega$ of full $\P$-measure such that these assertions hold for all $\omega \in \hat{\Omega}$. 
    
    To see that $A$ is $\P$-a.s.~invariant, we use the invariance of $\mu$ to find the identity
    \begin{equation}\label{eq:invariance_A_and_mu}
        A(\theta \omega) = \supp(\mu_{\theta\omega}) = \varphi_\omega( \supp(\mu_\omega))= \varphi_\omega(A(\omega)),
    \end{equation}
    for all $\omega \in \hat{\Omega}$. This shows that $A$ is $\P$-a.s.~invariant and, in summary, an invariant compact random set.  The fact that $A(\omega)$ is $\F_{-\infty}^0$-measurable follows directly from $\mu$ being a Markov measure.
    
    The main tool to prove that $A$ is a weak point attractor is \cite[Proposition 2.20]{Scheutzow17}, which implies that $A$ is a weak point attractor for points of the set
    \begin{equation}\label{eq:def_E_0}
        E_0:=\{ x : \lim_{n\to \infty} P^n(x,\cdot)\rightarrow \pi \}.
    \end{equation}
    If the associated Markov chain is aperiodic it is well known, e.g.~\cite[Theorem 2.27]{benaim2022markov}, that $E_0=\X$.
   If the Markov chain is periodic with period $p>1$, we divide the state space $\X$ into its $p$ periodic classes $W_1,\hdots, W_p$ such that $\varphi_\omega(W_i)\subset W_{i+i}$ $\P$-a.s.~where we identify $W_{p+1}$ with $W_1$, see \cite[Lemma 6.7.1]{durrett2019probability}. Let $\overline{\Omega} \subset \Omega$ be a $\theta$-invariant set of full $\P$-measure such that $\varphi_\omega(W_i)\subset W_{i+i}$ for all $\omega \in \overline{\Omega}$. A straightforward calculation using the Birkhoff ergodic theorem, cf.~\eqref{eq:Markov_chain_Birkhoff}, yields $\pi(W_i)=\frac{1}{p}$ for each $W_i$. By the definition of $\mu_\omega$ in (\ref{eq:correspondence_limit}), we get that $\mu_\omega(W_i)=\frac{1}{p}$ for all $\omega \in \overline{\Omega}$.
    
    On each set $W_i$, we consider the time-$p$ RDS over the probability space $(\overline{\Omega}, \F, \P)$, defined via the time-$p$ map $\varphi^{p}_\omega:W_i\rightarrow W_i$ and the time-$p$ shift $\theta_{p}$. Define the rescaled measures of $\mu$ restricted to $W_i$ by
    \begin{equation*}
        \mu_\omega^i:= \frac{1}{\mu_\omega(W_i)}\mu_\omega |_{W_i} = p \: \mu_\omega |_{W_i} .
    \end{equation*}
    This measure defines an invariant Markov measure of the time-$p$ RDS on $W_i$ since 
    \begin{equation*}
        (\varphi^p_\omega)_*(\mu_\omega|_{W_i}) = \mu_{\theta_p\omega}|_{W_i},
    \end{equation*}
    which follows immediately from the invariance of $\mu$. Therefore, the sets defined by 
    \begin{equation*}
        A^i(\omega):= \supp(\mu^i_\omega) = \supp(\mu_\omega) \cap W_i,
    \end{equation*}
    are invariant compact random sets that satisfy
    \begin{equation*}
        A^i(\omega)=A(\omega) \cap W_i.
    \end{equation*}
    The associated time-$p$ Markov chain on $W_i$ is aperiodic, irreducible and positive recurrent. Hence, the set $E_0^i$, defined analogously to~\eqref{eq:def_E_0}, equals $W_i$, and $A^i(\omega)$ is a weak point attractor of the time-$p$ RDS on $W_i$.
    For any point $x\in W_i$, we find 
    \begin{equation*}
        \lim_{n \to \infty} \P(\varphi^{pn}_\omega( x)\in A^i(\omega)) = 1.
    \end{equation*}
    Since $A^i(\omega)\subset A(\omega)$ and $\varphi^{N}_\omega(x)\in A(\omega)$ implies $\varphi^{n}_\omega(x)\in A(\omega)$ for any $n\geq N$, we find 
    \begin{equation*}
        \lim_{n \to \infty} \P(\varphi^{n}_\omega(x)\in A(\omega)) = 1.
    \end{equation*}
    This holds for any $x\in W_i$ and any class $W_i$. We conclude that $A$ is a weak point attractor.
\end{proof}
Step (II) has already been stated in \cite[Theorem 9]{EngelOlicon24}, for which we present a simplified proof.
\begin{proposition}[II]\label{prop:forward_point}
    $A$ is a forward point attractor.
\end{proposition}
\begin{proof}
    If $\varphi^N_\omega(x)\in A(\theta_N \omega)$, then also $\varphi^n_\omega(x)\in A(\theta_n \omega)$ for any $n\geq N$. Using continuity from below, we find
    \begin{align*}
        \P\big(\exists n\in \N: \varphi^n_\omega( x)\in A(\theta_n \omega)\big) &= \lim_{n\to \infty} \P\big(\varphi^n_\omega(x)\in A(\theta_n \omega)\big)=1.
    \end{align*}
    The second equality uses that $A$ is a weak point attractor. This shows that $A$ is a forward point attractor.
\end{proof}

\begin{proposition}[III]
    $A$ is a pullback point attractor.
\end{proposition}
\begin{proof}
 By the correspondence theorem, the pullback-limit of $\pi$ converges to $\mu_\omega$ for $\P$-a.e.~$\omega \in \Omega$, i.e.
    \begin{equation}\label{eq:weak_convergence}
        \lim_{n\to \infty} [(\varphi^{n}_{\theta_{-{n}}\omega})_*\pi] (B) = \mu_\omega(B), 
    \end{equation}
    for any $B\subset \X$. 
    By equality~(\ref{eq:A_supp_mu}), we have $\mu_{\omega} (A(\omega)) =1$  for $\P$-a.e.~$\omega$. Hence, equation~\eqref{eq:weak_convergence} yields 
    \begin{equation} \label{eq:corresp1}
        \lim_{n\to \infty} [(\varphi^{n}_{\theta_{-{n}}\omega})_*\pi] (A(\omega)) =1\,, \quad \text{ for $\P$-a.e.~$\omega$.}
    \end{equation}
    
    Assume now for contradiction that $A$ is not a pullback point attractor. Then there are an $x \in \X$ and a set $\Omega^x \subset \Omega$ with $\mathbb P (\Omega^x) > 0$ such that for every $\omega \in \Omega^x$ we have an increasing sequence $n_k \to \infty$ such that for each $k\in \N$
    \begin{equation*}
        \varphi^{n_k}_{\theta_{-{n_k}} \omega}(x) \notin A(\omega).
    \end{equation*}
This implies that $[(\varphi^{n_k}_{\theta_{-{n_k}} \omega})_*\pi](A(\omega))\leq 1 - \pi(x) < 1$ for any $\omega \in \Omega^x$ and $k \in \mathbb N$. This is a contradiction to~\eqref{eq:corresp1} since $\Omega^x$ has positive measure.
\end{proof}


\begin{proposition}[IV]\label{prop:point_attraction_equal_set}
    Pullback point attraction is equivalent to pullback attraction for finite deterministic sets and forward point attraction is equivalent to forward attraction of 
    compact random sets.
\end{proposition}
\begin{proof}
    Let $A$ be a pullback point attractor. We show that $A$ is a pullback attractor for finite deterministic sets.

    Let $B=\{x_1, \hdots, x_m \}\subset \X$ be a finite set. Since $A$ is a pullback point attractor, for $\P$-a.e.~$\omega \in \Omega$ there are numbers $N_1(\omega), \hdots, N_m(\omega)$ such that for all $n\geq N_i(\omega)$, we find $\varphi^n_{\theta_{-n}\omega} (x_i) \subset A(\omega)$. Let $N(\omega)=\max\{N_1(\omega), \hdots, N_m(\omega)\}$, implying $\varphi^n_{\theta_{-n}\omega} (B) \subset A(\omega)$ for all $n\geq N(\omega)$. 
    We conclude that
    \begin{equation*}
        \P(\exists N\in \N : \forall n\geq N: \varphi^n_{\theta_{-n}\omega} (B) \subset A(\omega)) = 1.
    \end{equation*}
    This shows that $A$ is a pullback attractor for finite deterministic sets. 
    
    Let now $A$ be a forward point attractor. By an analogous argument as above, $A$ is a forward attractor for finite deterministic sets.  
    Let $B$ be a compact random set and denote the set of finite subsets of $\X$ by $\mathcal{P}_c(\X)$, which is countable. 
    We compute
    \begin{align*}
        &\P(\exists N\in \N : \forall n\geq N: \varphi^n_{\omega} (B(\omega)) \subset A(\theta_n \omega))\\
        =\sum_{K \in \mathcal{P}_c(\X)} &\P(\exists N\in \N : \forall n\geq N: \varphi^n_{\omega} (K) \subset A(\theta_n \omega) \land B(\omega) = K)\\
        =\sum_{K \in \mathcal{P}_c(\X)} &\P(B(\omega) = K) = 1,
    \end{align*}
    where we used that $A$ is a forward attractor for finite deterministic sets. This shows that $A$ is a forward attractor for compact random sets.
\end{proof}

Step (V) is split into two propositions, the first one proving $(b)$ and the second proving $(c)$ and $(d)$.
\begin{proposition}[V]
    Let $A'$ be an invariant compact random set that satisfies some notion of attraction from Definition \ref{def:attractors}. Then $A'(\omega)=A(\omega)$ for $\P$-a.e.~$\omega \in \Omega$. 
\end{proposition}
\begin{proof}
    Weak point attraction is the weakest notion of attraction from Definition \ref{def:attractors}, and, hence, $A'$ is a weak point attractor. Applying Proposition \ref{prop:forward_point} and Proposition \ref{prop:point_attraction_equal_set}, to $A'$ instead of $A$, we conclude that $A'$ is a forward attractor for compact random sets. Hence, both $A$ and $A'$ are forward attractors for compact random sets, and consequently they attract each other. This means that for $\P$-a.e~$\omega \in \Omega$ there is an $N\in \N$ such that for all $n\geq N$
    \begin{equation*}
        A(\theta_n \omega) = \varphi_\omega^n(A(\omega)) \subset A'(\theta_n\omega), \qquad \text{and} \qquad A'(\theta_n \omega) = \varphi_\omega^n(A'(\omega)) \subset A(\theta_n \omega).
    \end{equation*}
    This shows that for $n\geq N$, we find $A'(\theta_n \omega) = A(\theta_n \omega)$. By ergodicity, we conclude that $A'(\omega) = A(\omega)$ for $\P$-a.e.~$\omega \in \Omega$.
\end{proof}

\begin{proposition}[V]
    There is a constant $\kappa\in \N$ such that $\abs{A(\omega)}= \kappa$ for $\P$-a.e.~$\omega\in \Omega$. The probability that a state $x\in \X$ lies in the attractor is given by
    \begin{equation*}
        \P(x\in A(\omega)) = \kappa \pi(x).
    \end{equation*}
\end{proposition}
\begin{proof}
    By \eqref{eq:sample_measures_uniform}, there is a constant $\kappa \in \N$ such that the sample measures of the unique invariant Markov measure $\mu$ are given by
    \begin{equation*}
        \mu_\omega = \frac{1}{\kappa} \sum_{a\in A(\omega)} \delta_a.
    \end{equation*}
    Hence, the attractor $A(\omega)$ consists of exactly $\kappa$ points for $\P$-a.e.~$\omega \in \Omega$. We compute
    \begin{equation*}
        \E[\mu_\omega(x)] = \frac{1}{\kappa} \P(x \in A(\omega)).
    \end{equation*}
    By the correspondence theorem, we have $\E[\mu_\omega] = \pi$, and the statement follows.
\end{proof}

\subsection{Size of the attractor}\label{sec:structure_of_attractor}
Let $(\theta, \varphi)$ be a fixed RDS representation of an irreducible and positive recurrent Markov chain with transition probabilities $P$ and stationary distribution $\pi$. In the previous section, we showed that the RDS has a unique attractor $A$ of size $\abs{A(\omega)}=\kappa$ for $\P$-a.e.~$\omega \in \Omega$. In the following, we characterize the size $\kappa$ of the attractor by analyzing the synchronizing behavior for sets of initial points.

We stress again that the RDS representation of a Markov chain is not unique. While the existence of a unique random attractor $A$ does not depend on the choice of the RDS representation (cf.~Theorem \ref{thm:equivalence}),
the properties of the attractor $A$, like its size, do (cf.~Example~\ref{ex:4_state_MC}).

We introduce the following relation on $\X$.
\begin{definition}\label{def:relation}
    We call two points $x,y\in \X$ \emph{insulated}, written $ x \parallel y$, if they do not $\omega$-synchronize for $\P$-a.e.~$\omega \in \Omega$, i.e.~
    \begin{equation*}
        \forall n\in \N : \P(\varphi^n_\omega(x)=\varphi^n_\omega(y))=0.
    \end{equation*}
    An \emph{insulated set} is a set $K\subset \X$ such that $x\parallel y$ for any $x\neq y\in K$. We denote the size of the largest insulated set by $\hat \kappa \in \N \cup \{\infty\}$. We call an insulated set $K\subset \X$ with $|K| = \hat \kappa$ a \emph{maximum insulated set}.
\end{definition}
\begin{remark}
\
\begin{enumerate}[(i)]
    \item There might be insulated sets $K\subset \X$ with $|K| < \hat \kappa$ such that adding any point $x$ to $K$ makes $K\cup \{x\}$ not insulated. Such a set might be called a \emph{maximal} insulated set, even though it is not a maximum insulated set. Maximal insulated sets will play no role in the rest of this work. Still, this distinction is important to avoid confusion in the terminology.
    \item Theorem \ref{thm:A_anti_clique} below shows that the $\P$-a.s.~size of the attractor $\kappa$ equals the size a maximum insulated set $\hat{\kappa}$. However, until this theorem is proven, we denote the size of a maximum insulated set by $\hat \kappa$ to distinguish it from the $\P$-a.s.~size of the attractor.
\end{enumerate}
    
\end{remark}
For an example of insulated sets of an RDS, see Example \ref{ex:4_state_MC}, in particular~\eqref{eq:Adams_RDS}. An important observation is that insulated sets $\P$-a.s.~remain insulated under the RDS. 
We state this fact without proof.
\begin{proposition}\label{prop:insulated_set_invariance}
    Let $K\subset \X$ be an insulated set. For any $n\in \N$ and $\P$-a.e.~$\omega\in \Omega$ the set $\varphi^n_\omega (K)$ is also an insulated set, with $|\varphi^n_\omega(K)|=|K|$.
\end{proposition}

The following theorem now establishes the random attractor as a maximum insulated set.
\begin{theorem}\label{thm:A_anti_clique}
    For $\P$-a.e.~$\omega \in \Omega$, the attractor $A(\omega)$ is a maximum insulated set and hence $\kappa = \hat\kappa$.
\end{theorem}
\begin{proof}
    By definition of a weak attractor, we have for any finite set $K$ that
    \begin{equation*}
    \lim_{n\to \infty} \P(\varphi^n_\omega(K) \subset A(\theta_n \omega)) =1.
    \end{equation*}
    If $K$ is a maximum insulated set, then, by Proposition \ref{prop:insulated_set_invariance}, for $\P$-a.e.~$\omega \in \Omega$  and any $n\in \N$ the set $\varphi^n_\omega(K)$ is a maximum insulated set as well. 
    Hence, using the measure-preserving property of $\theta$, the set $A(\omega)$ contains a maximum insulated set $\P$-a.s.~and, therefore, $|A(\omega)|\geq \hat \kappa$.
    
    In the proof of Proposition \ref{prop:existence_weak_attractor}, we already showed that $A(\omega)$ is $\P$-a.s.~of constant finite size $\kappa$. Assume now for contradiction that $\kappa >\hat \kappa$ for $\P$-a.e.~$\omega \in \Omega$. therefore, for each such $\omega$, the set $A(\omega)$ contains at least two points that are not insulated. Hence, with positive probability, the size of $A$ decreases in finite time, contradicting the $\P$-a.s.~fixed size of $A$. This shows that $\kappa= \hat{\kappa}$ and hence $A(\omega)$ is a maximum insulated set $\P$-a.s.
\end{proof} 

\begin{remark}\label{rem:kappa_sync}
\
\begin{enumerate}
    \item [(i)] Theorem \ref{thm:A_anti_clique} shows that the $\P$-a.s.~size $\kappa$ of the attractor equals $\hat \kappa$, the size of a maximum insulated set. Hence, we use $\kappa$ from now on to denote this quantity. Example \ref{ex:4_state_MC} illustrates that the value of $\kappa$ may differ for different RDS representations of the same Markov chain. The independent RDS always satisfies $\kappa= p$, where $p$ is the period of the Markov chain.
    \item[(ii)] The constant $\kappa$ determines whether a given RDS is synchronizing or not. If $\kappa >1$, the RDS cannot be synchronizing since points in a maximum insulated set $\P$-a.s.~never synchronize. If $\kappa = 1$, then any two points $x,y\in \X$ eventually hit the attractor $A(\omega)$ which consists of a single point. At that time $x$ and $y$ have synchronized. Hence, an RDS is synchronizing if and only if $\kappa = 1$. 
\end{enumerate}
\end{remark}

As mentioned before, and as demonstrated in Example \ref{ex:4_state_MC}, the $\P$-a.s.~size of the attractor $\kappa$ depends on the RDS representation of a given Markov chain. We provide restrictions on which values of $\kappa$ can be realized for a given Markov chain.

Recall that two points $x,y\in \X$ are said to $\omega$-synchronize if there is an $n\in\N$ such that $\varphi_\omega^n(x) = \varphi_\omega^n(y)$ (cf.~Definition \ref{def:omega_sync}). Given $\omega \in \Omega$, define the equivalence relation
\begin{equation*}
    x \sim_\omega y \quad \iff \quad \text{$x$ and $y$ $\omega$-synchronize}.
\end{equation*}
We denote the equivalence classes by $[x]_\omega$. We call these $\omega$-dependent equivalence classes synchronizing classes instead. Note that $x \parallel y$ is equivalent to $\P(x \sim_\omega y)=0$. In general, we can have $ 0 < \P(x \sim_\omega y) < 1$, cf.~Example \ref{ex:4_state_MC}.
\begin{proposition}\label{prop:coalescing_sets_uniform}
    For $\P$-a.e.~$\omega \in \Omega$, there are exactly $\kappa$ synchronizing sets. Each synchronizing set $[x]_\omega$ satisfies $\pi([x]_\omega) = \frac{1}{\kappa}$.
\end{proposition}
\begin{proof}
    Since the attractor $A$ is a forward point attractor, each $x\in\X$ $\omega$-synchronizes with one of the elements of $A(\omega)$, i.e. $x \in [a]_\omega$ for some $a\in A(\omega)$, for $\P$-a.e.~$\omega \in \Omega$. Since $A(\omega)$ consists of $\kappa$ elements for $\P$-a.e.~$\omega \in \Omega$, there are at most $\kappa$ synchronizing sets. On the other hand, insulated points $\P$-a.s.~lie in different synchronizing sets. Hence, for $\P$-a.e.~$\omega \in \Omega$ there are exactly $\kappa$ synchronizing sets.

    To prove the second statement, fix $\Tilde\omega \in \Omega$ such that there are exactly $\kappa$ synchronizing sets and each maximum insulated set has exactly one element inside each synchronizing set. This is the case for $\P$-a.e.~$\Tilde\omega$. Let $[x]_{\Tilde\omega}$ be one of the synchronizing sets. 
    By the correspondence theorem, we find $\pi([x]_{\Tilde\omega}) = \E[\mu_\omega([x]_{\Tilde\omega})]$, where the expectation is taken over $\omega$. The sample measures $\mu_{\omega}$ are a uniform distribution over the sets $A(\tilde\omega)$, cf.~\eqref{eq:sample_measures_uniform}, which yields
    \begin{equation*}
        \pi([x]_{\Tilde\omega}) = \E[\mu_\omega([x]_{\Tilde\omega})] = \int_{\Omega}\frac{1}{\kappa}\sum_{a\in A(\omega)} \delta_a([x]_{\Tilde\omega})\; \P(d\omega).
    \end{equation*}
    For $\P$-a.e.~$\omega \in \Omega$, the set $A(\omega)$ is a maximum insulated set. Hence, exactly one of the summands is $1$ while the others are $0$. We conclude $\pi([x]_{\Tilde\omega})=\frac{1}{\kappa}$. This holds for $\P$-a.e.~$\Tilde\omega\in \Omega$, which completes the proof.

\end{proof}

Proposition \ref{prop:coalescing_sets_uniform} imposes significant restrictions on the value of $\kappa$. This is relevant since the distribution $\pi$ depends solely on the associated Markov chain and not on the RDS representation.
\begin{corollary}\label{crl:restriction_kappa}
     If it is not possible to partition $\X$ into $k$ sets $W_i$ such that $\pi(W_i) = \frac{1}{k}$, then $\kappa \neq k$ for any RDS representation of the Markov chain.
\end{corollary}

\subsection{Examples}\label{sec:examples}
We provide an example of a Markov chain on a finite state space that corresponds to two different random dynamical systems with attractors of different cardinality, respectively.
\begin{example}\label{ex:4_state_MC} \textbf{(RDS representations with attractors of different sizes)}

    Let $\X = \{a, b, c, d\}$. Consider the Markov chain with transition probabilities $P$ as indicated by the following diagram.
    \begin{equation}\label{eq:Markov_chain_picture}
        \begin{tikzpicture}[node distance=1.5cm, ->, , baseline=(current  bounding  box.center)]
            \node(A)[]{$a$};
            \node(B)[right of=A]{$b$};
            \node(C)[below of=B]{$c$};
            \node(D)[below of=A]{$d$};
            
            \path (A) edge node[above] { $\frac{1}{2}$} (B);
            \path (B) edge node[right] {$ \frac{1}{2}$} (C);
            \path (C) edge node[below] {$ \frac{1}{2}$} (D);
            \path (D) edge node[left] {$ \frac{1}{2}$} (A);
            
            \path (A) edge [out=150,in=120,looseness=8] node[left] {$ \frac{1}{2}$} (A);
            \path (B) edge [out=60,in=30,looseness=8] node[right] {$ \frac{1}{2}$} (B);
            \path (D) edge [out=240,in=210,looseness=8] node[left] {$ \frac{1}{2}$} (D);
            \path (C) edge [out=330,in=300,looseness=8] node[right] {$ \frac{1}{2}$} (C);
        \end{tikzpicture}
    \end{equation}
    We present two possible RDS representations of this Markov chain. 

    We first consider the independent RDS $(\theta, \varphi)$, cf.~Section \ref{sec:canonical_constructions}. In each timestep any point may either take a step clockwise, or remain in its position. Since in the independent RDS all points move independently from one another, there are 16 possible maps $\varphi_\omega :\X \to\X$ that are taken with positive probability. Since all transition probabilities are $\frac{1}{2}$, each of these 16 possible maps occurs with probability $\frac{1}{16}$.
    It is not hard to verify that with probability $1$ all points eventually coalesce. Hence, the unique attractor $A(\omega)$ consists $\P$-a.s.~of a single point. 

    The second RDS representation we present is constructed over a shift space, like described in Section \ref{sec:canonical_constructions}. The map $\varphi_\omega:\X \to \X$ admits one of the following two configurations, each with probability $\frac{1}{2}$.
    \begin{equation}\label{eq:Adams_RDS}
         \begin{tikzpicture}[node distance=1.5cm, ->, baseline=(current  bounding  box.center)]
            \node(A1)[]{$a$};
            \node(B1)[right of=A1]{$b$};
            \node(C1)[below of=B1]{$c$};
            \node(D1)[below of=A1]{$d$};
            \node(f1)[right of=A1, xshift=-0.75cm, yshift=-0.75cm]{$f_1$};
            
            \path (A1) edge (B1);
            \path (C1) edge (D1);
            
            \path (B1) edge [out=60,in=30,looseness=8] (B1);
            \path (D1) edge [out=240,in=210,looseness=8] (D1);

            \node(A2)[right of=B1, xshift=2cm]{$a$};
            \node(B2)[right of=A2]{$b$};
            \node(C2)[below of=B2]{$c$};
            \node(D2)[below of=A2]{$d$};
            \node(f2)[right of=A2, xshift=-0.75cm, yshift=-0.75cm]{$f_2$};
            
            \path (B2) edge (C2);
            \path (D2) edge (A2);
            
            \path (A2) edge [out=150,in=120,looseness=8] (A2);
            \path (C2) edge [out=330,in=300,looseness=8] (C2);
        \end{tikzpicture}
    \end{equation}
    The RDS $(\theta, \varphi)$ over the shift space constructed by $\P_0(f=f_1)=\P_0(f=f_2) = \frac{1}{2}$ represents the Markov chain in~\eqref{eq:Markov_chain_picture}. Note that points diagonal to one another never coalesce. After one time step, all points lie on one of the two diagonals $\{a, c\}$ or $\{b, d\}$. Given $\omega=(\hdots, \omega_{-1}, \omega_0, \omega_1, \hdots) \in \Omega$ with $\omega_i\in \{f_1, f_2\}$, define
    \begin{equation*}
        A(\omega) = \begin{cases}
            \{a, c\}, \quad & \textup{if }\omega_{-1} = f_2,\\
            \{b, d\}, \quad & \textup{if }\omega_{-1} = f_1.
        \end{cases}
    \end{equation*}
    By construction, $A$ is a compact random invariant set. Additionally, $A$ absorbs any initial point of $\X$ within one time step. Hence, $A$ is the unique attractor which is $\P$-a.s.~of size 2.
    \hfill $\Diamond$
\end{example}

Theorem \ref{thm:equivalence} shows that for every RDS representing an irreducible, positive recurrent Markov chain there is an attractor $A$ that satisfies all notions of attraction from Definition \ref{def:attractors} except for pullback attraction for compact random sets. We illustrate via the following example that, in general, pullback attraction for compact random sets is not satisfied.
\begin{example}\textbf{(No pullback attraction for compact random sets)}\label{ex:no_pullback}

    Consider a random walk on $\X=\N$. The transition probabilities are given by $P(x, x+1) = \frac{1}{4}$, $P(x, x-1) = \frac{3}{4}$ for $x\geq 1$ and $P(0, 1) = \frac{1}{4}$, $P(0, 0)= \frac{3}{4}$. This Markov chain is irreducible and positive recurrent.
    
    We consider the RDS for which in each time step either all points move up or all points move down (except for the $0$, which stays at $0$). The map $\varphi_\omega$ admits one of the following two configurations with probability $\frac{1}{4}$ and $\frac{3}{4}$ respectively:
    \begin{equation*}
        f_1(x) = x+1, \qquad f_{-1}(x) = \max\{x-1, 0\}.
    \end{equation*}
    Let $(\Omega, \F, \P)$ be the probability space over which the RDS $(\theta, \varphi)$ is defined. 
    By Theorem~\ref{thm:equivalence}, the RDS admits a unique attractor $A$ which is $\P$-a.s.~a singleton and we write $A(\omega) = \{a(\omega)\}.$ We construct a compact random set $B(\omega) = \{b(\omega)\}$ such that $A$ does not attract $B$ in the pullback sense. Note that for the chosen RDS the distance of two points can decrease by at most one in each step. Hence, if $|b(\theta_{-n}\omega) - a(\theta_{-n} \omega) | > n $, then $\varphi_{\theta_{-n} \omega}^n b(\theta_{-n} \omega) \neq a(\omega)$. We construct the point $b(\omega)$ in dependence of $a(\omega)$ by
    \begin{equation*}
        b(\omega) \coloneqq a(\omega) + g(a(\omega)),
    \end{equation*}
    where $g:\N \to \N$ is an increasing function. Constructing $b(\omega)$ in this way implies that $b$, and thereby the set $B$, are $\F_{-\infty}^0$-measurable.

    For $k,m \in \N$ consider the probabilities
    \begin{equation*}
        \P\big(a(\theta_{-n} \omega) \leq k, \; \forall n \in \{0, \hdots, m\}\big).
    \end{equation*}
    For fixed $k\in \N$ and $m \to \infty$, this probability goes to $0$, because the sequence $a(\theta_{-n} \omega)$ is $\P$-a.s.~unbounded. Therefore, it is possible to construct a strictly increasing $g:\N \to \N$ such that for all $k\in \N$
    \begin{equation*}
        \P\big(a(\theta_{-n} \omega) \leq k, \; \forall n \in \{0, \hdots, g(k)\}\big) \leq \frac{1}{k^2}.
    \end{equation*}
    In particular the probabilities are summable over $k$. By the Borel--Cantelli theorem, with probability 1, only finitely many of these events occur. This means that for $\P$-a.e.~$\omega \in \Omega$, there is a $K\in \N$ such that for all $k\geq K$ we have
    \begin{equation*}
        \exists n \in \{0, \hdots, g(k)\}: \; a(\theta_{-n}\omega) > k.
    \end{equation*}
    This gives a sequence of $n_k$, where $k\geq K$, such that $a(\theta_{-n_k}\omega) > k$. By construction, the $n_k$ satisfy $n_k \leq g(k)$. Using that $g$ is strictly increasing, we obtain
    \begin{equation*}
        b(\theta_{-n_k} \omega) = a(\theta_{-n_k} \omega) + g(a(\theta_{-n_k} \omega)) > a(\theta_{-n_k} \omega) + g(k) \geq a(\theta_{-n_k} \omega) + n_k.
    \end{equation*}
    By the argument above, that points can only get closer to each other by at most one in each step, we conclude
    \begin{equation}\label{eq:pullback_fail}
        \varphi_{\theta_{-n_k}}^{n_k} b(\theta_{-n_k} \omega) \neq a(\omega).
    \end{equation}
    This holds for all $n_k$, for $k\geq K$. A priori, the sequence $(n_k)_{k\in \N}$ does not need to be increasing, nor do the values have to be distinct. However, $ a(\theta_{-n_k} \omega) \geq k$ such that the sequence $n_k$ has to be diverging. In particular, there are infinitely many $n_k \in \N$ for which \eqref{eq:pullback_fail} holds. This shows that $A$ does not attract $B$ in the pullback sense.
    \hfill $\Diamond$
\end{example}

\section{Expected hitting times} \label{sec:time_until_attraction}
Let $(\theta, \varphi)$ be a fixed RDS representation of an irreducible and positive recurrent Markov chain with transition probabilities $P$ and stationary distribution $\pi$. Theorem \ref{thm:equivalence} states that there is a unique attractor $A$. If the attractor is $\P$-a.s.~a singleton, i.e.~$\kappa=1$, the RDS is synchronizing, see Remark \ref{rem:kappa_sync}. For synchronizing RDS, the synchronization time of two initial conditions $x,y \in \X$ is defined as
\begin{equation}\label{eq:def_T_Delta}
    T_\Delta(\omega, x,y ) = \min\{n\in \N \mid \varphi_\omega^n(x) = \varphi_\omega^n(y) \},
\end{equation}
and is $\P$-a.s.~finite for any $x,y\in \X$. For general RDS, which are not necessarily synchronizing, define the random hitting time of the attractor for an initial condition $x\in \X$,
\begin{equation}\label{eq:def_attraction_time}
    T_A(\omega, x) := \min \{n\in \N\: | \: \varphi^n_\omega(x) \in A(\theta_n \omega) \}.
\end{equation}
Since $A$ is a forward attractor, the random time $T_A(\omega, x)$ is $\P$-a.s.~finite.

In this section, we study whether the synchronization time and the time until attraction have finite expectation. The main results are the following. Theorem \ref{thm:synchronization_finite_mean} states that for synchronizing RDS the expectation $\E[T_\Delta(\omega, x,y)]$ is always finite for all $x,y\in \X$, while Theorem \ref{thm:finite_mean_iff_erg_2} shows that the time until attraction has finite expectation for all initial conditions given that the associated Markov chain is ergodic of degree 2. In Example \ref{ex:ergodic_deg_2}, we show that without the assumption of ergodicity of degree 2, the expected hitting time of the attractor can be infinite. 

A Markov chain is ergodic of degree 2 if for some $x\in \X$ we have
\begin{equation*}
    \E_x[\tau^2_x] < \infty,
\end{equation*}
where $\tau_x$ is the first hitting time of $x$ as given in \eqref{eq:MC_hitting_time}. 
The following proposition lists equivalent characterizations of ergodicity of degree 2 which we will use in the following.
For alternative characterizations of ergodicity of degree 2 and higher degrees, we refer to \cite{mao2003algebraic}.
\begin{proposition}\label{prop:ergodic_degree_2}
    Consider an irreducible, positive recurrent Markov chain with stationary distribution $\pi$. The following are equivalent:
    \begin{enumerate}[(i)]
        \item The Markov chain is ergodic of degree 2.
        \item $\E_x[\tau_y^2] < \infty$, for all $x,y\in \X$.
        \item $\E_\pi[\tau_y] = \sum_{x\in \X} \pi(x) \E_x[\tau_y] < \infty$, for all $y\in \X$.
    \end{enumerate}
\end{proposition}
\begin{proof}
    The equivalence of $(i)$ and $(ii)$ is shown in~\cite[Remark 3.1]{mao2003algebraic}. We show the equivalence of $(i)$ and $(iii)$. Due to~\cite[Proposition 2.10]{benaim2022markov}, we have the identity
    \begin{equation*}
        \E_\pi[\tau_y] =\frac{1}{2} \pi(y)(\E_y[\tau^2_y] + \E_y[\tau_y]). 
    \end{equation*}
    Since $\E_y[\tau_y] < \infty$ for positive recurrent Markov chains, we find $\E_\pi[\tau_y]<\infty$ if and only if $E_y[\tau_y^2] < \infty$. 
\end{proof}
We now state the main results of this section.

\begin{theorem}\label{thm:synchronization_finite_mean}
    Let $(\theta, \varphi)$ be a synchronizing Markov RDS representing an irreducible and positive recurrent Markov chain. For any $x,y\in \X$, the expected synchronization time is finite, i.e.~
    \begin{equation*}
        \E[T_\Delta(\omega, x,y)] < \infty.
    \end{equation*}
\end{theorem}

\begin{theorem}\label{thm:finite_mean_iff_erg_2}
    Assume that the associated Markov chain is ergodic of degree 2. Then the expected hitting time of the attractor $\E[T_A(\omega, x)]$ is finite for all $x\in \X$, and we have 
    \begin{equation}\label{eq:time_attraction_finite}
       \sum_{x \in \X} \pi(x) \E[T_A(\omega, x)] < \infty.
    \end{equation}
\end{theorem}

\begin{remark}
    Ergodicity of degree 2 is a property of the Markov chain itself while the expected time until attraction depends on the RDS representation.
    However, even for Markov chains which are ergodic of degree 2, the time until attraction may still be arbitrarily high depending on the RDS representation. In other words, there is no upper bound on the quantity in \eqref{eq:time_attraction_finite} which does not depend on the RDS representations. This is demonstrated in Example \ref{ex:slow_sync}.
\end{remark}

Theorem \ref{thm:finite_mean_iff_erg_2} shows that ergodicity of degree 2 is a sufficient condition for the attraction time to have finite mean. We believe that ergodicity of degree 2 is also a necessary condition. Up to this point, we have not found a proof for this conjecture.
\begin{conjecture}\label{con:finite_one_finite_all}
    If the associated Markov chain is not ergodic of degree 2, then $\E[T_A(\omega,x)] = \infty$ for all $x\in \X$.
\end{conjecture}

Both Theorem \ref{thm:synchronization_finite_mean} and Theorem \ref{thm:finite_mean_iff_erg_2} are consequences of a more general statement about the coalescence and insulation of (random) initial conditions. We generalize the definition of the synchronization time to not necessarily synchronizing RDS
\begin{equation} \label{eq:general_T_Delta}
    T_\Delta(\omega, x,y ) := \min\{n\in \N \mid \varphi_\omega^n(x) = \varphi_\omega^n(y) \: \lor \: \varphi_\omega^n(x) \parallel  \varphi_\omega^n(y)\}.
\end{equation}
Note that for synchronizing RDS, no two points are insulated such that this definition coincides with \eqref{eq:def_T_Delta} and there is no conflict of notation.
\begin{proposition}\label{prop:decision_finite_mean}
    Let $(\theta, \varphi)$ be a Markov RDS representing an irreducible and positive recurrent Markov chain. Let $T_\Delta(\omega, x,y)$ be as defined in \eqref{eq:general_T_Delta}. The following holds:
    \begin{enumerate}[(I)]
        \item For any $x,y\in \X$, we find
        \begin{equation*}
            \E[T_\Delta(\omega, x,y)] < \infty.
        \end{equation*}
        \item If the associated Markov chain is ergodic of degree 2, we have
        \begin{equation*}
            \sum_{x,y \in \X} \pi(x) \pi(y) \E[T_\Delta(\omega, x,y)] < \infty.
        \end{equation*}
    \end{enumerate}
\end{proposition}
Theorem \ref{thm:synchronization_finite_mean} is an immediate consequence of Proposition \ref{prop:decision_finite_mean} (I). For a proof of Proposition \ref{prop:decision_finite_mean}, and how Theorem \ref{thm:finite_mean_iff_erg_2} follows from (II), see Section \ref{sec:proofs_attraction_time}.

\subsection{Examples}
The first example is a Markov chain that is not ergodic of degree 2 for which the expected hitting time of the attractor is infinite.
\begin{example}\label{ex:ergodic_deg_2}
    \textbf{(Infinite expected attraction time)} Consider the Markov chain on $\X = \N_{> 0}$ with transition probabilities $P(n, 1)= p_n$ and $P(n, n+1) = 1 - p_n$, where $p_n \in (0,1)$. This Markov chain is clearly irreducible. We claim that the probabilities $p_n$ can be chosen in a way such that $\E_1[\tau_1] < \infty$, but $\E_1[\tau_1^2] = \infty$. Indeed, given any probability distribution $\varrho$ on $\N_{>0}$ with full support, we can define $p_n \coloneqq \P(\varrho = n \mid \varrho \geq n)$ which results in $\tau_1\sim \varrho$. 
    
    We consider the independent RDS $(\theta, \varphi)$, cf.~Section \ref{sec:canonical_constructions}, which has $\kappa=1$, i.e.~the attractor is a singleton $A(\omega) = \{a(\omega)\}$. Note that the only way that two trajectories coalesce is if both of them jump to $1$ in the same time step. Hence, for $x\in \X$, unless $a(\omega) = x$, we have
    \begin{equation*}
        T_A(\omega, x) \geq \min \{n\in \N \mid \varphi_\omega^n(a(\omega)) = 1 \}.
    \end{equation*}
    The random singleton $a(\omega)$ is distributed according to $\pi$, cf.~\eqref{eq:prob_x_in_A}. We conclude
    \begin{equation*}
        \E[T_A(\omega, x)] \geq \sum_{a \in \X} \pi(a) \E_a[\tau_1] - \E_x[\tau_1].
    \end{equation*}
    By Proposition \ref{prop:ergodic_degree_2}, the right-hand side is infinite, since we constructed the Markov chain to be not ergodic of degree 2. This shows that for any $x\in \X$, the expected time of hitting the attractor is infinite.
    \hfill $\Diamond$
\end{example}

The next example provides RDS representations of a Markov chain, which is ergodic of degree 2, with arbitrarily high expected times until attraction.
\begin{example}\label{ex:slow_sync}
    \textbf{(Arbitrarily high mean attraction times)} Consider the Markov chain on two states $\X=\{a,b\}$ with transition probabilities $P(a,a) = P(a,b) = P(b,a) = P(b,b) = \frac{1}{2}$. Given $\epsilon\in(0,1)$, we construct an RDS representation over a shift space, cf.~Section \ref{sec:canonical_constructions}, such that
    \begin{align*}
        \P(\varphi_\omega(a,b) &= (a,b)) = \frac{1-\epsilon}{2},\\
        \P(\varphi_\omega(a,b) &= (b,a)) = \frac{1-\epsilon}{2},\\
        \P(\varphi_\omega(a,b) &= (a,a)) = \frac{\epsilon}{2},\\
        \P(\varphi_\omega(a,b) &= (b,b)) = \frac{\epsilon}{2}.
    \end{align*}
    At any time, the probability for the two states to coalesce is exactly $\epsilon$.
    Hence, maximum insulated sets have size $\kappa = 1$. The expected time to hit the attractor is geometrically distributed. We have
    \begin{equation*}
        \E[T_A(\omega, a)] = \frac{1}{2}\E[T_A(\omega, a) \:|\: A(\omega) = b] = \frac{1}{2\epsilon}.
    \end{equation*}
    As $\epsilon \to 0$ this expectation becomes arbitrarily large.

    The same technique works for any Markov chain that has an RDS representation with $\kappa >1$. Given an RDS representation with $\kappa>1$, fix $\epsilon>0$ and let all states move independently from one another with probability $\epsilon$. The resulting RDS has $\kappa = 1$, but the expected time until states which were previously insulated coalesce is at least $\frac{1}{\epsilon}$. An interesting open question is whether there is an upper bound on the expected hitting time of the attractor, in particular of the quantity in \eqref{eq:time_attraction_finite}, which is uniform over all RDS representations, in the special case that the Markov chain has no RDS representation with $\kappa >1$. Note that Corollary~\ref{crl:restriction_kappa} provides a sufficient condition for a Markov chain to be not representable by an RDS with $\kappa>1$.  \hfill $\Diamond$
\end{example}

\subsection{Proofs}\label{sec:proofs_attraction_time}
This section contains the proofs of Proposition \ref{prop:decision_finite_mean} and Theorem \ref{thm:finite_mean_iff_erg_2}.

\subsubsection{Proof of Proposition~\ref{prop:decision_finite_mean}}
    Let $P$ be the transition probabilities of the 
    associated Markov chain and denote its stationary distribution by $\pi$.

    Consider the Markov chain on $\X^2$ which models the 2-point motion of the RDS $(\theta, \varphi)$. It is defined by the transition probabilities
    \begin{equation}\label{eq:def_Q}
        Q((x,y), (x',y')) := \P (\varphi_\omega(x,y) = (x', y')).
    \end{equation}
    To distinguish probabilities and expected values of the Markov chain on $\X^2$ defined by $Q$ from those of the Markov chain on $\X$ defined by $P$, we write $\P^Q$ and $\E^Q$.
    Consider the set
    \begin{equation*}
        \Delta := \{(x,y)\in \X^2 \mid x=y \: \lor \: x\parallel y\}.
    \end{equation*}
     Since $Q$ models the 2-point motion of the RDS, we find 
    \begin{equation}\label{eq:exp_T_Delta_Q}
        \E[T_{\Delta}(\omega, x^*, y^*)] = \E_{(x^*, y^*)}^Q[\tau_{\Delta}],
    \end{equation}
    for any $x^*, y^*\in \X$, where $\tau_{\Delta}$ is the hitting time of the set $\Delta$. 
    

    \begin{proof}[Proof of Proposition~\ref{prop:decision_finite_mean} (I)]
    Fix $x^*, y^* \in \X$. Without loss of generality, we assume that the associated Markov chain is aperiodic. If there was a period $p$, we distinguish two cases. If $x^*$ and $y^*$ lie in different periodicity classes, then $x^* \parallel y^*$ and $\tau_\Delta(\omega, x^*,y^*)=0$. If $x^*$ and $y^*$ lie in the same periodicity class, consider the $p$-step Markov chain, which is aperiodic.
    
    Consider the Markov chain on $\X^2$ which models the 2-point motion of the RDS $(\theta, \varphi)$ with transition probabilities $Q$, as defined in \eqref{eq:def_Q}. We denote a trajectory of the Markov chain defined by $Q$ by $(\mathbf{X}_n)_{n\in \N}$, where $\mathbf{X}_n = (\mathbf{X}_n^1, \mathbf{X}_n^2)$. Note that both components evolve according to the transition probabilities $P$ of the 1-point Markov chain. Let $C \subset \X$ be a finite set such that $\pi(C)>\frac{1}{2}$. Using the Birkhoff ergodic theorem for Markov chains, cf.~Theorem~\ref{thm:MC_ergodic}, we compute
    \begin{nalign}\label{eq:Birkhoff_C_argument}
        \lim_{N\to \infty} \frac{1}{N} \sum_{n=1}^{N} \1_{C^2}(\mathbf{X}_n)
        \geq&\lim_{N\to \infty} \frac{1}{N} \sum_{n=1}^{N} \left(1 - \1_{\X \setminus C} (\mathbf{X}_n^1) - \1_{\X \setminus C} (\mathbf{X}_n^2) \right)\\
        =&1 - 2 (1-\pi(C)) > 0.
    \end{nalign}
    Hence, any trajectory of $Q$ spends a significant part of its time in the finite set $C^2$. In particular, every irreducibility class $I$ of $Q$ must intersect $C^2$. This implies that there are only finitely many irreducibility classes. We denote the union of all irreducibility classes by $\mathcal{I}$. Since all pairs of initial conditions eventually either coalesce or become insulated, we find $\mathcal{I} \subset \Delta$. 

    By aperiodicity of $P$, for every $c\in \X$ and $z\in \X$ there is an $N\in \N$ such that for all $n\geq N$ we find $P^n(c,z ) >0$. This implies that there is an $N \in \N$ such that for each $c\in C$, we find $P^N(c,x^*)>0$ and $P^N(c,y^*) > 0$. We fix such an $N$ and define 
    \begin{equation}\label{eq:def_altered_R}
        R((x,y), (x', y')) := \begin{cases}
            Q^N((x,y), (x', y')), \quad & (x,y) \notin \mathcal I,\\
            P^N(x, x') P^N(y, y'), \quad & (x,y) \in \mathcal I.
        \end{cases}
    \end{equation}
    This defines the transition probabilities of another Markov chain on $\X^2$. When outside of the set $\mathcal{I}$, the transition probabilities are equivalent to taking $N$ steps with respect to $Q$. Inside $\mathcal{I}$, the two components move independently from one another for $N$ steps. Thus, the set $\mathcal{I}$ is not an irreducibility class under $R$.
    
    We claim that $R$ has a unique irreducibility class $J$, which contains $(x^*, y^*)$. Note that the two individual components of trajectories governed by $R$ obey the transition probabilities $P^N$. Let $(\mathbf{Y}_n)_{n\in \N} = (\mathbf{Y}_n^1, \mathbf{Y}_n^2)$ be a trajectory of the Markov chain defined by $R$. By the same argument as \eqref{eq:Birkhoff_C_argument}, $\mathbf{Y}_n$ spends a significant portion of its time in $C^2$. In particular, any initial condition $(x,y) \in \X^2$ can reach $C^2$ under $R$. By construction, any point $(c_1, c_2) \in C^2$ has a positive probability to land in $(x^*, y^*)$ after one step of $R$. This shows that the Markov chain defined by $R$ has a unique irreducibility class $J$ that contains $(x^*, y^*)$. Since $(\mathbf{Y}_n)_{n\in \N}$ spends a significant part of its time in $C^2$, there is at least one $\mathbf c \in C^2$ such that
    \begin{equation}\label{eq:c_is recurrent}
        \lim_{N\to \infty} \frac{1}{N} \sum_{n=1}^{N} \1_\mathbf{c}(\mathbf{Y}_n) > 0.
    \end{equation}
    A classic result from Markov chain theory, e.g.~\cite[Proposition 2.1]{benaim2022markov}, shows that \eqref{eq:c_is recurrent} implies that $R$ is positive recurrent on $J$.
    
    Outside of $\mathcal{I}$, the transition probabilities $R$ coincide with $Q^N$. Since $\mathcal{I} \subset \Delta$, we find
    \begin{equation}\label{eq:Expectation_Q_R}
        \E_{(x^*, y^*)}^Q[\tau_\Delta] \leq \E_{(x^*, y^*)}^Q[\tau_\mathcal{I}] \leq N \E_{(x^*, y^*)}^R[\tau_\mathcal{I}].
    \end{equation}
    We showed that $R$ is positive recurrent on its unique irreducibility class $J$, $(x^*, y^*) \in J$ and $J\cap \mathcal{I} \neq \emptyset$. This implies that the right-most term of \eqref{eq:Expectation_Q_R} is finite. We conclude that $\E_{(x^*, y^*)}^Q[\tau_\Delta]$ is finite. Since $\E_{(x^*, y^*)}^Q[\tau_\Delta] = \E[\tau_\Delta(\omega, x^*, y^*)]$, cf.~\eqref{eq:exp_T_Delta_Q}, this completes the proof.
    \end{proof}
    
    \begin{proof}[Proof of Proposition~\ref{prop:decision_finite_mean} (II)]
    Let $C\subset \X$ be a finite set with $\pi(C) > \frac{5}{6}$. Define the random time until two initial conditions reach $C^2$.
    \begin{equation*}
        T_{C^2}(\omega, x,y) = \min \{n \in \N \mid \varphi_\omega^n(x,y) \in C^2 \}.
    \end{equation*}
    The first part of the proposition states that the expected time $\E[T_\Delta(\omega, c_1, c_2)]$ is finite for each fixed pair $(c_1,c_2) \in C^2$. Since $C^2$ is a finite set, we have
    \begin{equation*}
        \max_{(c_1, c_2) \in C^2} \E[T_\Delta(\omega, c_1, c_2)] < \infty.
    \end{equation*}
    Hence, it suffices to show that randomly chosen initial conditions $(x,y)$ enter $C^2$ in a time of finite mean since
    \begin{equation}\label{eq:split_sum}
        \sum_{x,y \in \X} \pi(x) \pi(y) \E[T_\Delta(\omega, x,y)] \leq \sum_{x,y \in \X} \pi(x) \pi(y) \E[T_{C^2}(\omega, x,y)] + \max_{(c_1, c_2) \in C^2} \E[T_\Delta(\omega, c_1, c_2)].
    \end{equation}

    Let $x\in \X$ and consider the trajectory $X_n = \varphi_\omega^n(x)$. This is a trajectory of the associated Markov chain with transition probabilities $P$. Let $c \in C$ and set $\gamma \coloneqq \frac{4}{5}\pi(c)$. Define
    \begin{equation*}
        N_n = \sum_{k=1}^n \1_{c}(X_k), \qquad  R_n = N_n-\gamma n.
    \end{equation*}
    By the Birkhoff theorem, cf.~Theorem \ref{thm:MC_ergodic}, we find $\lim_{n\to \infty} \frac{1}{n} R_n = \pi(c) -\gamma > 0$. Consider the return times to the point $c$ which are recursively defined by 
    \begin{align*}
        T_0 &= \min \{n \in \N \mid X_n = c \}, \\
        T_{k+1} &= \min \{n> T_k \mid X_n = c \}.
    \end{align*}
    The process $S_k = R_{T_k}$ is a random walk on $\R$ with $S_0 = 1-\gamma T_0$. The increments $\xi_k := S_k - S_{k-1}$ are i.i.d.~and given by $\xi_k = 1-\gamma(T_k - T_{k-1})$. Thus, the increments are distributed like
    \begin{equation*}
        \xi \sim 1 -\gamma \tau_{c},
    \end{equation*}
    where $\tau_c$ is the random return time from $c$ to itself. We conclude $\E[\xi] = \frac{1}{5}$ and $\E[\abs{\xi}^2] < \infty$ from ergodicity of degree 2. 
    
    Consider the absolute value of the minimum of the random walk, relative to its starting point
    \begin{equation*}
        m = -\min_{n\in \N} (S_n - S_0) \geq 0.
    \end{equation*}
    By \cite[Chapter X.2, Theorem 2.1]{Asmussen1998}, $m$ is almost surely finite, and we have $\E[m] < \infty$ given that $\E[\abs{\xi}^2] < \infty$.

    The sequence $R_n$ is decreasing except at times $T_k$ when it increases by $1-\gamma$. At these times its value is given by $R_{T_k} = S_k$. We conclude $R_n \geq S_0 -m-1$ for all $n\in \N$. This implies
    \begin{equation*}
        \frac{1}{n}N_n = \frac{1}{n} R_n + \gamma \geq \frac{S_0-m-1}{n} + \gamma.
    \end{equation*}
    Define the random time
    \begin{equation*}
        T_{\frac{3}{5}c}(\omega, x) \coloneqq 4\frac{m+1 - S_0}{\gamma}.
    \end{equation*}
    This time is defined such that for all $N\geq T_{\frac{3}{5}c}$, we find
    \begin{equation}\label{eq:3_5_c}
        \frac{1}{N} \sum_{n=1}^N \1_c(X_n) = \frac{1}{N} N_N \geq \frac{3}{4} \gamma = \frac{3}{5} \pi(c).
    \end{equation}
    The quantities $S_0$ and $m$ are random variables. However, note that the expected value of $m$ does not depend on the initial condition $x$. That is because $m$ only considers the minimum of $S_n$ after the trajectory has first hit $c$. The expected value of $S_0$ is given by $\E[S_0] = 1-\gamma \E_x[\tau_c]$. If the initial condition $x$ is randomly distributed with respect to $\pi$, we conclude
    \begin{equation}\label{eq:3_5_random_x}
        \sum_{x\in \X} \pi(x) \E[T_{\frac{3}{5}c}(\omega, x)] \leq 4 \frac{\E[m] + \gamma \E_\pi [\tau_c]}{\gamma} <  \infty.
    \end{equation}
    The term $\E_\pi [\tau_c]$ is finite because the Markov chain is assumed to be ergodic of degree 2, cf.~Proposition \ref{prop:ergodic_degree_2}.

    There are only finitely many $c\in C$. Define
    \begin{equation*}
        T_{\frac{3}{5}C}(\omega, x) = \max_{c\in C} T_{\frac{3}{5}c}(\omega, x).
    \end{equation*}
    By \eqref{eq:3_5_random_x} the expectation of each $T_{\frac{3}{5}c}(\omega, x)$ is finite when summed over $\pi(x)$. Hence, the same holds for $T_{\frac{3}{5}C}(\omega, x)$, i.e.
    \begin{equation}\label{eq:T_C_finite}
        \sum_{x\in \X} \pi(x) \E[T_{\frac{3}{5}C}(\omega, x)] < \infty.
    \end{equation}
    From \eqref{eq:3_5_c}, we obtain that for all $N\geq T_{\frac{3}{5}C}(\omega, x) $,
    \begin{equation*}
        \frac{1}{N}\sum_{n=1}^N \1_C(X_n) \geq \frac{3}{5} \pi(C) > \frac{1}{2},
    \end{equation*}
    where we used the assumption $\pi(C) > \frac{5}{6}$. In conclusion, we have shown that for an initial condition, randomly chosen with respect to $\pi$, there is a time of finite mean after which the Birkhoff sums of the indicator function $\1_C$ are always greater than $\frac{1}{2}$. 
    
    Consider two initial conditions $x,y\in \X$. For $N \geq \max\{T_{\frac{3}{5}C}(\omega, x), T_{\frac{3}{5}C}(\omega, y)\}$, both trajectories have spent more than half of the time inside the set $C$. Hence, there is at least one $n\leq N$ for which $X_n = \varphi_\omega^n(x) \in C$ and $Y_n=\varphi_\omega^n(y) \in C$. We find
    \begin{equation*}
        T_{C^2}(\omega, x,y) \leq \max\{T_{\frac{3}{5}C}(\omega, x), T_{\frac{3}{5}C}(\omega, y)\}.
    \end{equation*}
    From \eqref{eq:T_C_finite} we conclude
    \begin{align*}
        \sum_{(x,y)\in \X^2} \pi(x) \pi(y) \E[T_{C^2}(\omega, x,y)] &\leq \sum_{(x,y)\in \X^2} \pi(x) \pi(y) \E[T_{\frac{3}{5}C}(\omega, x) + T_{\frac{3}{5}C}(\omega, y)] \\
        &= 2 \sum_{x\in \X} \pi(x) \E[T_{\frac{3}{5}C}(\omega, x)] < \infty.
    \end{align*}
    This shows that the right-hand side of \eqref{eq:split_sum} is finite, and completes the proof.
    \end{proof}

\subsubsection{Proof of Theorem \ref{thm:finite_mean_iff_erg_2}}
    Assume that the associated Markov chain is ergodic of degree 2 and let $x \in \X$. Note that $\varphi_\omega^n(x) \in A(\theta_n \omega)$ if and only if for all $a \in A(\omega)$, we have either $\varphi_\omega^n(x) = \varphi_\omega^n(a)$ or $\varphi_\omega^n(x) \parallel \varphi_\omega^n(a)$. More precisely, $x$ has coalesced with exactly one of the $a$'s and is insulated from the others. We conclude
    \begin{equation*}
        T_A(\omega, x) = \max_{a \in A(\omega)} T_\Delta(\omega, x, a),
    \end{equation*}
    where $T_\Delta$ is as defined in \eqref{eq:general_T_Delta}.
    This enables us to estimate the expectation of $T_A(\omega, x)$ via considering the times $T_\Delta(\omega, x,y)$ for pairs $(x,y)\in \X^2$:
    \begin{align*}
        \E[T_A(\omega, x)] &\leq \E\left[ \sum_{a\in A(\omega)} T_\Delta(\omega, x, a) \right] =  \E\left[ \sum_{y\in \X} \1_{A(\omega)}(y) T_\Delta(\omega, x, y) \right]\\
        &= \sum_{y\in \X} \E[ \1_{A(\omega)}(y) T_\Delta(\omega, x, y)].
    \end{align*}
    Note that all quantities are non-negative such that the exchange of an infinite sum and the expectation is justified.
    Since $\1_{A(\omega)}$ is $\F_{-\infty}^0$-measurable, cf.~Theorem \ref{thm:equivalence}, and $T_\Delta(\omega, x, y)$ is $\F_0^\infty$-measurable, we obtain 
    \begin{nalign}\label{eq:bound_tau_A}
        \E[T_A(\omega, x)] &\leq \sum_{y\in \X} \E[ \1_{A(\omega)}(y)] \E[T_\Delta(\omega, x, y)] \\
        &= \sum_{y\in \X} \kappa \pi(y) \E[T_\Delta(\omega, x, y)],
    \end{nalign}
    recalling that $\P(y\in A(\omega)) = \kappa \pi(y)$, cf.~\eqref{eq:prob_x_in_A}. In view of \eqref{eq:time_attraction_finite}, we have
    \begin{equation}\label{eq:bound_tau_A_2}
        \sum_{x\in \X} \pi(x) \E[T_A(\omega, x)] \leq \kappa \sum_{x,y\in \X} \pi(x) \pi(y) \E[T_\Delta(\omega, x, y)].
    \end{equation}
    By Proposition \ref{prop:decision_finite_mean} (II), the right-hand side is finite, which completes the proof. \hfill $\square$

\section{Discussion}\label{sec:discussion}


In this work, we have explored the fundamental structure of random attractors for discrete-time random dynamical systems on discrete state spaces.
In Section \ref{sec:equivalence}, we have shown the existence of a unique random attractor $A$, that is a forward attractor for compact random sets and a pullback attractor for finite deterministic sets. However, $A$ is typically not a pullback attractor for compact random sets. It is common to define pullback attraction for a class of 
suitable compact random sets, e.g.~tempered random sets \cite[Definition 14.3]{RasmussenKloeden11}. Investigating which classes of compact random sets are pullback attracted by $A$ --- with or without assuming ergodicity of degree 2 --- is a topic for potential future research.
Our investigations in Section \ref{sec:structure_of_attractor} have structured the state space via insulated and synchronizing sets, showing that the $\P$-a.s.~size $\kappa$ of the attractor coincides with the size of any maximum insulated set. 

Section \ref{sec:time_until_attraction} has characterized expected hitting times: we have shown that for any RDS representation of a Markov chain that is ergodic of degree 2, any initial condition hits the attractor in a time of finite mean. The converse statement has been left as a conjecture.
Furthermore, under ergodicity of degree 2, the expected time to hit the attractor, when the initial condition is distributed according to $\pi$, is finite as well (see formula~\eqref{eq:time_attraction_finite}).
An interesting open question is which RDS representation minimizes this quantity, i.e.~for which RDS representation of a given Markov chain with randomly distributed initial conditions the attractor is reached within the shortest time span. 
Identifying this RDS would be particularly useful for sampling from the stationary distribution $\pi$, using the method of the coupling from the past algorithm \cite{propp1996exact}.

\section*{Acknowledgements}
This research has been partially supported by Deutsche Forschungsgemeinschaft (DFG) through grant CRC 1114 ``Scaling Cascades in Complex Systems'', Project Number 235221301, Project A02 ``Multiscale data and asymptotic model assimilation for atmospheric flows '' (G.O.-M.) and Project A08 ``Characterization and prediction of quasi-stationary atmospheric states'' (R.C. and M.E.), and through Germany's Excellence Strategy -- The Berlin Mathematics Research Center MATH+ EXC-2046/1, project 390685689, in particular subprojects AA1-8 (M.E. and G.O.-M.) and AA1-18 (M.E.). M.E. additionally thanks the Einstein Foundation and the NWO (VI.Vidi.233.133) for supporting his research.
R.C. thanks Adam Schweitzer for engaging discussions, and providing valuable examples, like \eqref{eq:Adams_RDS} in Example \ref{ex:4_state_MC}. The authors thank Francois Baccelli, Dennis Chemnitz and Tobias Hurth for helpful discussions and contributing valuable ideas. The authors also thank the anonymous reviewers for their outstanding suggestions, leading to significant simplifications of the proofs of Section \ref{sec:time_until_attraction}, and for providing the idea for Example \ref{ex:no_pullback}.

\end{document}